\documentclass{amsart} 
\usepackage{epsf} 
\usepackage{amssymb,latexsym, amsmath, amscd, array} 
\def\hepsffile{\leavevmode\epsffile} 
\usepackage{graphicx} 

\swapnumbers 
\numberwithin{equation}{section} 

\theoremstyle{plain} 
\newtheorem{thm}{Theorem}[section]

\newtheorem{lem}[thm]{Lemma}

\newcommand\theoref{Theorem~\ref} 
\newcommand\lemref{Lemma~\ref}

\newcommand\defref{Definition~\ref}

\theoremstyle{definition} 
\newtheorem{defin}[thm]{Definition}

\newtheorem{rem}[thm]{Remark} 
\newtheorem{ex}[thm]{Example}

\def\id{\protect\operatorname{id}} 
\def\Im{\protect\operatorname{Im}}

\def\rk{\protect\operatorname{rank}} 
 
\def\Res{\protect\operatorname{Res}}

\def\win{\protect\operatorname{awin}} 
\def\pass{\protect\operatorname{pas}} 
\def\aawin{\widetilde{\win}}
\def\aawin{\protect\operatorname{AWIN}} 
\def\clk{\protect\operatorname{lk}} 
\def\alk{\protect\operatorname{alk}} 
\def\swin{\protect\operatorname{win}}

\def\eps{\varepsilon} 
\def\gf{\varphi}

\def\C{{\mathbb C}} 
\def\Z{{\mathbb Z}} 
\def\Q{{\mathbb Q}} 
\def\R{{\mathbb R}} 
 
\def\bA{{\pmb A}}  
\def\bB{{\pmb B}}

\def\1{\hbox{\rm\rlap {1}\hskip.03in{\rom I}}} 
\def\Bbbone{{\rm1\mathchoice{\kern-0.25em}{\kern-0.25em} 
{\kern-0.2em}{\kern-0.2em}I}} 

\def\pp{\medskip{\parindent 0pt \it Proof.\ }} 
\def\wt{\widetilde} 
 
\def\m{\smallskip} 
\def\ov{\overline} 

\def\A{{\mathcal A(M, \mathcal N)}} 
\def\B{{\mathcal B(M, \mathcal N)}} 

\begin{document} 
\date{November 7, 2006} 
\leftline{ } 
\centerline{ } 

\title[On generalized winding numbers] 
{On generalized winding numbers} 
\author[V.~Chernov (Tchernov) and Yu.~Rudyak]{Vladimir V. 
Chernov 
(Tchernov) and Yuli B. Rudyak} 
\address{V. Chernov, Department of Mathematics, 
6188 Bradley Hall, Dartmouth College, Hanover NH 03755, 
USA} 
\email{Vladimir.Chernov@dartmouth.edu} 
\address{Yu. Rudyak, Department of Mathematics, University 
of Florida, 358 
Little Hall, Gainesville, FL 32611-8105, USA} 
\email{rudyak@math.ufl.edu} 



\begin{abstract} 
Let $M^m$ be an oriented manifold, let $N^{m-1}$ be an oriented closed manifold, 
and let $p$ be a point in $M^m$. For a smooth map $f:N^{m-1} \to M^m, p\notin 
\Im f,$ we 
introduce an invariant $\win_p(f)$
that can be regarded as a generalization of 
the classical winding number of a planar curve around a point. We show 
that $\win_p$ estimates from below the number of times a wave 
front  on 
$M$ passed 
through a given point $p\in M$ between two moments of time.
Invariant $\win_p$  allows us to formulate the analogue of 
the complex analysis Cauchy integral formula for meromorphic 
functions on complex surfaces of genus bigger than 
one.
\end{abstract}

\maketitle 

\section*{Introduction}
Gauss linking number is a link homotopy invariant of a pair 
$(\phi_1(N_1^{n_1}), \phi_2(N_2^{n_2}))$ of disjoint linked closed oriented  
submanifolds  in an oriented manifold $M^m$ of dimension $m=n_1+n_2+1.$ 
The linking number $\clk$ is defined via basic homology theory as the 
intersection number of a singular chain whose boundary is $\phi_2(N_2)$ with 
$\phi_1(N_1).$

If $\phi_{2*}([N_2])\neq 0\in H_*(M),$ then $\phi_2(N_2)$ is not a boundary of 
any singular chain. If $\phi_{2*}([N_2])=0$ but $\phi_{1*}([N_1])\neq 0\in 
H_*(M),$ then the intersection 
number in the definition above does depend on the relative homology class in 
$H_*(M, \phi_2(N_2))$ realized by the singular chain.
Thus $\clk$ is not well defined unless the link components are zero homologous.
(Similar homology theory methods allow one to define the linking number in the 
case
where $\phi_{i*}([N_i])\in H_*(M)$ are of finite order or 
$\phi_{i*}([N_i])=0\in H_*(M, \partial M), i=1,2,$~see 
Kaiser~\cite{Kaiserbook}.)

The winding number $\swin_pf$ of a loop
$f: S^1 \to \R^2$ around a point $p\notin \Im f$ measures how many times $f$ 
turns around $p$. It is defined as the intersection number of a path $P$ 
connecting $p$ to infinity with a loop $f$.
Clearly the winding number is just the linking number between $f: S^1\to 
S^2=\R^2\cup\{\infty\}$ and the map of $S^0$ sending one point of $S^0$ to $p$ 
and the other to $\infty$
.

The winding number is a classical invariant. It is a part of the complex 
analysis Cauchy integral formula and it has many applications in topology. The 
result of Whitney~\cite{Whitney} expresses the rotation number of a planar curve 
through the winding numbers of the curve around points in the regions of the 
curve complement. 
Many  formulas involving the winding number for invariants of planar curves, 
fronts, and knot diagrams were obtained in the works of Polyak~\cite{Polyak}, 
Shumakovich~\cite{ShumakovichStrangeness},~\cite{ShumakovichVassiliev}, 
Turaev~\cite{Turaevshadow}, Viro~\cite{Viro} and the first 
author~\cite{Tchernov}. Generalizations of winding numbers to hypersurfaces in 
$\R^m$ are essential in the works of Goryunov~\cite{Goryunov}, 
Mikhalkin and Polyak~\cite{MikhalkinandPolyak}.

Since the winding number is a particular case of the linking number, standard 
homology theory methods give its generalization to the case of a (possibly 
singular) smooth hypersurface $f:N^{m-1}\to M^m$, $p\not \in \Im f,$ provided 
that $f_*([N])=0\in H_{m-1}(M)$ and that $\partial M\neq \emptyset$ or that $M$ 
is the interior of a manifold with boundary.
(A boundary component plays the role of infinity where one places the second 
point of $S^0.$ The condition $f_*([N])=0$ is needed so  that the 
intersection number of the path $P$ and the hypersurface $f(N)$ does not depend 
on the choice of the path $P$.) {\it We denote the winding number defined this 
way by $\swin_p(f).$\/}

In our work~\cite{CR} we constructed the ``affine linking invariant'' $\alk$ of 
a pair of linked singular oriented closed submanifolds $\bigl(\phi_1(N_1^{n_1}), 
\phi_2(N_2^{n_2})\bigr)$ of $M^{n_1+n_2+1}.$ Our $\alk$ is a link homotopy invariant. It is a 
generalization of the linking number and it is well defined for all 
$\phi_{1*}([N_1]), \phi_{2*}([N_2])\in H_*(M).$ The group where the $\alk$ 
invariant takes values depends on the homotopy classes of the maps $\phi_1, 
\phi_2$ and sometimes is hard to compute.

In this paper we use ideas, similar to the ones we used in~\cite{CR} to define 
$\alk,$ to construct the {\it affine winding number\/} $\win_p(f)$ that is a 
generalization of the winding number to a vast collection of oriented $M$ and 
closed oriented $N.$ We do not require that $f_*([N])=0\in H_*(M)$ and that 
$\partial M\neq \varnothing$ (or that $M$ is the interior of a manifold with 
boundary). Thus in these cases the invariant $\swin_p(f),$ that is a particular 
case of the linking number, is not defined.

Since one of the linked manifolds is the one-point-space, many of the technical 
difficulties we dealt with in~\cite{CR} do not arise. In particular, the group 
where the affine winding number takes values is either $\Z$ or a quotient group 
of $\Z.$ Also  the operation $\mu$ on the bordism groups, that we introduced 
in~\cite{CR} to define $\alk,$  reduces to the standard intersection pairing. 
(The operation $\mu$ is quite interesting on its own. It gives rise to a Poisson 
bracket on the bordism group of 
garlands,~\cite{ChernovRudyakGarlands},~\cite{ChernovPoisson} that is related to 
Goldman-Turaev~\cite{Goldman},~\cite{Turaevskein} and 
Andersen-Mattes-Reshetikhin~\cite{AMR1},~\cite{AMR2} algebras.)

The generalized affine winding numbers we construct in this paper have affine 
nature, i.e. only the difference of affine winding 
numbers of two homotopic maps is well-defined. Equivalently, 
for a fixed distinguished $\eps:N\to M$ we can define the 
affine winding number $\win_p(f)$ of $f: N \to M$ around $p$ provided that $f$ 
and $\eps$ are homotopic. 
In the classical case of $M=\R^2$ and $N=S^1$ such a  distinguished map is, in a 
sense, fixed 
implicitly and it is a map into a point in $\R^2\setminus p$.
The precise definitions are given in 
Section \ref{winding}.  

\m
In Section \ref{apps} we consider two applications of our theory. First, we 
formulate the generalization of the complex analysis Cauchy integral formula to 
the 
case of meromorphic functions on complex surfaces of genus bigger than 
one. 
Second, we consider a propagation of a wave 
front in $M$ and estimate the number of passages of the front through a  given 
point between two time moments 
$t_1$ and $t_2$. In many cases we can estimate this number just from the shapes 
of the 
front at time-moments $t_1$ and $t_2,$ without any 
knowledge of the propagation process, topology of $M$, etc. 

\m
In Section \ref{versiya} we generalize the invariant $\win_p(f)$  to the case 
where the point $p\in M$ is not fixed and  somehow moves  in $M.$

\section{Affine winding numbers}\label{winding} 
We work in the $C^{\infty}$-category and the word {\it smooth\/} means 
$C^{\infty}.$
In this paper $N^{m-1}$ and $M^m$ are oriented connected smooth manifolds of 
dimensions
$m-1$ and $m$ respectively, and we assume $m\geq 2.$ The manifold $N$ is assumed 
to be closed. $I$ is the interval $[0,1].$

In this section we fix a point $p\in M$ and a connected component $\mathcal N$ 
of the space $C^{\infty}(N, M)$ of smooth maps $N\to M.$ Put $\Sigma\subset 
\mathcal N$ to be the {\it discriminant\/} that consists of all maps $f\in 
\mathcal N$ such that $p\in \Im f.$ We do not include into $\Sigma$ the maps $f$ 
that are singular in the common sense but do not satisfy $p\in \Im f.$ 

The affine winding number $\win_p(f)$ that we define in this section is a 
locally constant function on $\mathcal N\setminus \Sigma.$ Equivalently, it is a 
function on $\pi_0(\mathcal N\setminus \Sigma).$ If $f_*([N])=0\in H_*(M)$ and 
$\partial M\neq \emptyset,$ so that $\swin_p(f)$ can be defined as the 
particular case of the linking number, then the functions $\swin_p(f), 
\win_p(f):\mathcal N\setminus \Sigma\to \Z$ are equal up to an additive 
constant.

\begin{defin}\label{good}
Let $p$ be a point in $M$.
We say that a smooth map $F: N\times I\to M$ is {\it good} if $p\notin 
F(N\times\{0,1\})$ and $p$ is a regular value of $F$. We also call such $F$ {\it 
a good homotopy\/} between
the maps  $F|_{N\times 0}, F|_{N\times 1}\in \mathcal N\setminus \Sigma.$
\end{defin}

\m If two smooth maps $N\to M$ are homotopic and their images do not contain 
$p$, then 
there exists a good homotopy between these maps. Moreover, the set of good 
homotopies is $C^0$-dense in the set of all homotopies. This 
can be proved via standard general position arguments.

\m Note that $F^{-1}(p)$ is a finite set for every good $F$. The 
standard orientation of $I$ yields an orientation of $N\times I.$ 
Every point in $F^{-1}(p)$ is equipped with a sign $\pm 1$ as follows. We put 
the {\it sign of the point\/} to be $+1$ if the restriction of $F$ to a small 
neighborhood of the point is orientation preserving, and we put the sign of the 
point to be $-1$ otherwise.

\begin{defin}\label{delwin}
For a good map $F:N\times I\to M,$ we define $\Delta_{\win_p}(F)\in \Z$ to be 
the  sum of the signs of the inverse images of 
$p$ under $F$.
\end{defin}

\begin{rem}[relation between $\Delta_{\win_p}(F)$ and the intersection 
number]\label{alter}
The number $\Delta_{\win_p}(F)$ can also be described as follows. Consider the 
maps 
$$
\Phi: N\times I \to M\times I, \, \Phi(x,t)=(F(x),t)
$$ 
and 
$$
P: I \to  M\times I, P(t)=(p,t).
$$ 
Then $\Delta_{\win_p}(F)$ is equal to the intersection 
number of $\Phi$ and $P$. The proof is straightforward.

Similarly $\Delta_{\win_p}(F)$ equals to the intersection number of $F$ and of 
the positively oriented point $p.$
\end{rem}

\m We regard the circle $S^1$ as the quotient space $I/\{0,1\}$ and denote by 
$\pi : I \to S^1$ the projection. Consider a good map $F:  N \times I \to M$  
and assume that $F(x,0)=F(x,1),$ for all $x\in N$. Then there exists the unique 
map 
$G: N\times S^1\to M$ with $F=G\circ (1_N\times\pi)$.
The following Lemma  follows immediately from Definition \ref{delwin} and the 
description of 
the degree as the sum of the signs of the inverse images. {\it Note that if the 
manifold $M$ is not closed, then the degree of $G$ is zero by definition.\/}

\begin{lem}\label{degree}
$\Delta_{\win_p}(F)$ equals to the degree of the map $G$.\qed
\end{lem}

\begin{defin}[$\A$]\label{A} 
We call a smooth map
$\mu:N\times S^1\rightarrow M$ {\it special} if $\mu\big |_{N\times s}\in 
\mathcal N,$  for some (and therefore for all) $s\in S^1$. We define the {\it 
indeterminacy subgroup} 
$\bA=\bA(M, \mathcal N)$ of $\Z$ to be the  subgroup of possible degrees of 
special 
maps $\mu:N\times S^1\rightarrow M$. Put $\A=\Z/\bA$ and denote by 
$q:\Z \to \A$ the quotient homomorphism. 
\end{defin} 

\begin{defin} \label{win}
Fix $\eps=\eps_{\mathcal N}\in\mathcal N\setminus \Sigma$ that should be thought 
of as a preferred map.
Given a map $f\in \mathcal N\setminus \Sigma$, we define the {\it affine winding 
number\/} 
$\win_p(f)=\win_{p, \eps}(f)\in \A$ by setting 
$$
\win_p(f)=q(\Delta_{\win_p}(F))\in \A,
$$ 
where $F:N\times I \to M$ is a good homotopy between $\eps=F|_{N\times 0}$ and 
$f=F|_{N\times 1}.$ 
\end{defin}

\m  Note that $\win_{p,\eps}(\eps)=0$.

\begin{thm}\label{welldefined}
The invariant $\win$ has the following properties:
\begin{description}
\item[1]
The number $\win_p(f)=\win_{p, \eps}(f)\in \A$ does not depend on the choice of 
the good homotopy $F$ in its definition. Thus $\win_p:\mathcal N\setminus \Sigma 
\to \A$ is well defined. 
\item[2] If $H$ is a good homotopy between $f_0=H|_{N\times 0}$ and 
$f_1=H|_{N\times 1},$
then $\win_p(f_1)-\win_p(f_0)=q(\Delta_{\win_p}(H)).$ In particular, 
$\win_p:\mathcal N\setminus \Sigma \to \A$ is a locally constant function.
\item[3] If we replace $\eps$ by some $\eps'\in \mathcal N\setminus \Sigma,$ 
then 
$\win_{p, \eps'}-\win_{p,\eps}:\mathcal N\setminus \Z\to \A$ is the constant 
function $\win_{p, \eps'}(\eps).$ Hence the functions $\win_{p, \eps'}, 
\win_{p,\eps}:\mathcal N\setminus \Z\to \A$ are equal up to an additive 
constant.
\end{description}
\end{thm}

\pp {\it Let us prove statement $1.$\/} Let $F_1, F_2:N\times I\to M$ be two 
good homotopies between $\eps$ and $f.$ We must show that 
$q(\Delta_{\win_p}(F_1))=q(\Delta_{\win_p}(F_2))\in \A.$
Define $\ov F_2:N\times I\to M$ by $\ov F_2(n, t)=F_2(n, 1-t).$ Clearly 
$\Delta_{\win_p}(\ov F_2)=-\Delta_{\win_p}(F_2).$ Thus it suffices to show that 
$q\bigl (\Delta_{\win_p}(F_1)+\Delta_{\win_p}(\ov F_2)\bigr)=0\in \A.$

A homotopy $F_1$ followed by $\ov F_2$ gives a homotopy $F:N\times I\to M$ 
defined via $F(n, t)=F_1(n, 2t),$ for $n\in N, t\in [0,\frac{1}{2}],$ and via 
$F(n, t)=\ov F_2(n, 2t-1),$ for $n\in N, t\in [\frac{1}{2}, 1].$ Perturbing $F$ 
slightly in $N\times (0,1),$ if needed, we can assume that it is smooth and 
hence that $F$ is a good homotopy between $F|_{N\times 0}=\eps$ and $F|_{N\times 
1}=\eps.$ Clearly, $\Delta_{\win_p}(F)= \Delta_{\win_p}(F_1)+\Delta_{\win_p}(\ov 
F_2).$
\lemref{degree} implies that $q(\Delta_{\win_p}(F))=q\bigl 
(\Delta_{\win_p}(F_1)+\Delta_{\win_p}(\ov F_2)\bigr)=0\in \A.$ Thus the function 
$\win_p:\mathcal N\setminus \Sigma\to \A$ is well defined. 

{\it To prove statement $2,$\/} choose a good homotopy $F:N\times I\to M$ from 
$F|_{N\times 0}=\eps$ to $F|_{N\times 1}=f_0.$ Similarly to above, $F$ followed 
by $H$ gives a good homotopy $J$ from $\eps$ to $f_1.$
Clearly $\Delta_{\win_p}(J)=\Delta_{\win_p}(F)+\Delta_{\win_p}(H).$ By 
definition of $\win_p,$ we have 
$\win_p(f_1)=q(\Delta_{\win_p}(J))$ and $\win_p(f_0)=q(\Delta_{\win_p}(F)).$
Thus 
$\win_p(f_1)-\win_p(f_0)=q\bigl(\Delta_{\win_p}(J)-\Delta_{\win_p}(F)\bigr)=
q(\Delta_{\win_p}(H)).$ 

If $f_0$ and $f_1$ belong to the same path connected component of $\mathcal 
N\setminus \Sigma,$ then we can find a good homotopy $H$ between them such that 
$H^{-1}(p)=\emptyset.$ Thus $\win_p(f_1)-\win_p(f_0)=q(0)=0$ and hence $\win_p$ 
is a locally constant function.

{\it Let us prove statement $3.$\/} Clearly $\A$ and $q:\Z\to \A$ do not depend 
on the choice of $\eps\in \mathcal N\setminus \Sigma.$ Choose
a good homotopy $F':N\times I\to M$ from $\eps'$ to $\eps.$ 
Choose $f\in \mathcal N\setminus \Sigma$ and a good homotopy $F$ from $\eps$ to 
$f.$ We have
$\win_{p, \eps}(f)=q(\Delta_{\win_{p}}(F)).$ 
Similarly to above, homotopy $F'$ followed by $F$ gives a good homotopy from 
$\eps'$ to $f.$ Counting the preimages of $p$ under this homotopy we get that 
$\win_{p, \eps'}(f)=q(\Delta_{\win_{p}}(F))+q(\Delta_{\win_{p}}(F')).$ Thus 
$\win_{p, \eps'}(f)-\win_{p,\eps}(f)=q(\Delta_{\win_{p}}(F'))=
\win_{p, \eps'}(\eps).$\qed

\begin{rem}[Affine nature of $\win_p$]\label{differ} 
By statement $3$ of Theorem~\ref{welldefined}, if we change the distinguished 
map $\eps\in \mathcal N\setminus \Sigma,$ then the function $\win_p=\win_{p, 
\eps}:\mathcal N\setminus \Sigma\to \A$ changes by an additive constant. Thus if 
we neglect the distinguished map $\eps,$ then the invariant $\win_p$ is 
well-defined up to an additive constant.
This is similar to the ambiguity in the choice of the origin in an affine space, 
and this shows the affine nature of our invariant $\win_p.$

Note also that statement $2$ of Theorem~\ref{welldefined} implies that for $f_0, 
f_1\in \mathcal N\setminus \Sigma$ the difference $\win_{p, \eps}(f_1)-\win_{p, 
\eps}(f_0)$ does not depend on the choice of the distinguished map $\eps.$
\end{rem}

\m It is useful to know when the indeterminacy subgroup $\bA$ is trivial, i.e. 
when
$\A=\Z$. For such spaces, similarly to the classical case of $M=\R^2$ and 
$N=S^1$, the affine winding numbers are indeed integer numbers rather than  
elements of $\A=\Z/\bA$. 

\begin{thm}\label{m=0} 
The equality $\A=\Z$ holds if and only if all the special mappings $N\times 
S^1\to M^m$ have zero degrees. 
In particular, $\A=\Z$ provided 
that at least one of the following conditions is satisfied: 
\begin{description} 
\item[0] the manifold $M^m$ is not closed.
\item [1] there exists $i\in \{1,\ldots, m\}$ such that $b_i(N)+b_{i-1}(N)<b_i 
(M)$, where $b_i$ is the $i$-th Betti number.
\item [2] the space $\mathcal N$ consists of  null-homotopic maps and 
there exists $i\in \{1, \ldots, m\}$ 
such that $b_{i-1}(N)<b_i(M)$. $($In particular, this condition holds 
if $b_1(M)>1$, since $N$ is connected.$)$ 
\item [3] $N=S^{m-1}$, the space $\mathcal N$ consists of null homotopic maps, 
and $M$ is not a rational homology sphere.
\item [4] a map from $\mathcal N$ induces the trivial 
homomorphism $\pi_1(N)\to \pi_1(M),$ the group $\pi_1(M)$ is infinite and 
it does not 
contain $\Z$ as a finite index subgroup.
\item [5] $M$ is a closed manifold that admits a complete Riemannian metric of 
negative sectional curvature.

\end{description} 
\end{thm} 

\pp Clearly $\A=\Z$ if and only all the special maps $\mu:N\times 
S^1\to 
M$ have zero degrees.
Thus it suffices to show that if any of the 
conditions $0-5$ is 
satisfied, then every special mapping $\mu :N\times S^1\to M$ 
has zero degree. 
We deal with each condition separately. 

{\bf Condition $0$.\/} The degree 
of a map from a closed $m$-dimensional manifold to a non-closed $m$-dimensional 
manifold is always zero. Thus every special $\mu:N\times S^1\to M$ has zero 
degree. 

{\em For this reason, while considering the cases of Conditions 
{\bf 1--5} we can 
and shall assume that $M$ is closed.\/} 
 
{\bf Condition $1$.\/} Consider a special map $\mu:N\times 
S^1\to M$ of  
degree $l \neq 0$. We must prove that $b_i(N)+b_{i-1}(N) 
\geq b_i(M)$ for all 
$i$. 

Let $\mu_!:H_*(M) \to H_*(N\times S^1)$ be the transfer 
map, see 
e.g. \cite[V.2.11]{Rudyak}. As it is well-known, $\mu_*\circ 
\mu_! (x)=lx$ for all $x\in H_*(M)$. In particular, 
$ 
\mu_*: H_*(N \times S^1;\Q) \to H_*(M;\Q) 
$ 
is an epimorphism. 
Hence $\rk H_i(N\times S^1)=b_i(N\times S^1)\geq \rk H_i(M) 
=b_i(M)$, for all 
$i=1,\dots, m$. By the K\"unneth formula we have 
$b_i(N\times S^1)=b_i(N)+b_{i-1}(N)$. Thus $b_i(N)+b_{i-1} 
(N)\geq b_i(M)$ 
for all $i$.

{\bf Condition $2$.\/} Consider a special map $\mu:N\times 
S^1\to M$ of  
degree $l \neq 0$. We must prove that $b_{i-1}(N)\geq b_i
(M) 
$ for all $i$. 

Similarly to the case of Condition {\bf 1}, we conclude 
that 
$ 
\mu_*: H_*(N \times S^1;\Q) \to H_*(M;\Q) 
$ 
is an epimorphism. 
Fix a point $*\in S^1$ and denote by $i: N \to N \times 
S^1, i(n)=(n,*)$ the 
inclusion. Consider a diagram of homology groups 
$$ 
\CD 
H_i(N; \Q) @>i_*>>  H_i(N\times 
S^1; \Q) @>\mu_*>> 
H_i(M;\Q)\\ 
@.  @VVp_*V @.\\ 
@.  H_i(N;\Q) @. 
\endCD 
$$ 
where $p:N\times S^1\to N$ is the projection. 

Since $\mathcal N$ consists of null-homotopic maps, $\mu_*i_*$ is the zero 
homomorphism. Because $p_*i_*=\id,$ we get that $i_*$ is injective and thus 
$\rk \ker \mu_*\geq b_i(N)$. Since $\mu_*$ is surjective, we get that $\rk \Im 
\mu_*=b_i(M).$ Using K\"unneth formula and elementary linear algebra we have
$b_i(N)+b_{i-1}(N)=b_i(N\times S^1)=\rk \ker \mu_* +\rk \Im \mu_*\geq 
b_i(N)+b_i(M).$
Thus we have $b_{i-1}(N)\geq b_i(M)$. 

{\bf Condition $3$.\/} Assume that $\dim M=m>2$. If $M$ is 
not a rational 
homology sphere, then there exists $i$ with $1\le i\le m-1$ 
and such that 
$b_i(M)>0$. Moreover, we can take $i>1$ since $b_1(M)=b_{m- 
1}(M)$ by the 
Poincar\'e duality. So 
$ 
0=b_{i-1}(S^{m-1})<b_i(M) 
$ and condition {\bf 2} holds. 

If $m=2$ and $M$ is not a rational homology sphere, then 
$M \ne S^2$ is a closed oriented surface. So $b_1(M)>1=b_0(N)$ 
and again 
condition {\bf 2} holds. 

{\bf Condition $4$.\/} Let $\mu:N\times S^1\to M$ be a special 
map of degree 
$l \neq 0.$ Consider the cover map $p:\widetilde 
M\to M$ such that 
$$ 
\Im \bigl (p_*:\pi_1(\widetilde M)\to \pi_1(M) \bigr) =\Im 
\bigl( \mu_*:\pi_1(N\times S^1) \to \pi_1(M)\bigr). 
$$ 
Let $\widetilde {\mu}:N\times S^1\to \widetilde M$ be a $p$-lifting of $\mu$, 
i.e. a
map such that 
$\mu=p \circ \widetilde {\mu}$. 
Clearly $\deg(\widetilde \mu)\deg(p)=\deg (\mu)=l\neq 0,$ where $\deg$ denotes 
the degree of a map. 
Hence the covering 
$p$ should be finite and 
$\Im  \mu_*$ 
is a finite index 
subgroup of $\pi_1(M)$. Since $\mu$ is special and a map from $\mathcal N$ 
induces the trivial homomorphism $\pi_1(N)\to \pi_1(M),$
the 
composition 
of the homomorphisms
$$ 
\CD 
\pi_1(N) = \pi_1(N \times \{*\})@>i_{1*}>>\pi_1(N \times S^1) @>\mu_*>> \pi_1(M) 
\endCD 
$$  
is the trivial homomorphism. 
Now using the equality 
$ 
\pi_1(N\times S^1)=\pi_1(N) \oplus \pi_1(S^1) 
$  and the fact that $\Im \mu_*$ is a finite index subgroup of $\pi_1(M),$ we 
conclude that 
$$ 
\CD 
\Im \Bigl ( \pi_1(S^1)=\pi_1(\{ *\} \times S^1) @>i_{2*}>> 
\pi_1(N\times 
S^1) @>\mu_*>> 
\pi_1(M)\Bigr ) 
\endCD 
$$ 
is a finite index subgroup of $\pi_1(M)$. 

If $\Im (\mu_* \circ i_{2*})$ is finite, then $\pi_1(M)$ is finite. If $\Im 
(\mu_* \circ i_{2*})$ is infinite, then it is isomorphic to $\Z$ and therefore 
$\pi_1(M)$ contains $\Z$ as a finite index subgroup.

{\bf Condition $5$.\/}  By the Hadamard Theorem~\cite{DoCarmo} $M$ is a 
$K(\pi_1(M),1)$-space. 
Take a special map $\mu:N\times S^1\to M$. 

First, assume that the map 
$$
\CD
h: S^1=\{*\}\times S^1 @>i_2>> N \times  S^1 @>\mu >> M
\endCD
$$ is homotopy trivial.   Since $M$ is a $K(\pi_1(M),1)$-space, standard 
obstruction theory arguments show that $\mu$ is homotopic to a map that passes 
through a projection onto $N^{m-1}$.  Thus $\deg \mu=0$.

Now assume that the 
map $h$
is homotopy non-trivial.  Clearly $\pi_1(N\times \{*\})$ 
commutes 
with $\pi_1(\{*\}\times S^1)$ in $\pi_1(N \times S^1)$.  The Preissman 
Theorem~\cite{DoCarmo} says that all nontrivial abelian subgroups of $\pi_1(M)$ 
are infinite cyclic. Therefore $\Im \mu_*(\pi_1(N 
\times S^1))$ 
is an infinite cyclic subgroup of $\pi_1(M)$ and the
homomorphism 
$\mu_*$ has the form $\pi_1(N\times S^1) \to \Z \to \pi_1(M)$. 
Furthermore, 
the  inclusion $\Z\subset \pi_1(M)$ can be induced by a map 
$\psi:  S^1 \to M$. Since $M$ is a $K(\pi_1(M), 1)$-space, standard obstruction 
theory arguments show that  
$\mu: S^1 \times N \to M$ is homotopic to a map 
$ 
\overline\mu: S^1\times N \rightarrow S^1 \xrightarrow{\psi} M. 
$ 
Thus, $\deg \mu=\deg \overline \mu=0$.
This completes the proof of \theoref{m=0}.
\qed 

\begin{rem} 
It is not always true that $\A=\Z.$
For example, if $M=T^2,$ $N=S^1$ and $\mathcal N$ consists of maps homotopic to 
the 
meridian, then $\A=0.$  The same is true for $M=S^2, N=S^1$.
\end{rem}
 
\m  
The following Theorem says that if the winding number $\swin_p$ can be defined 
as the particular case of the classical linking number, then the invariants 
$\swin_p, \win_p:\mathcal N\setminus \Sigma\to \Z$ are equal up to an additive 
constant.   
Note that using conditions $1-5$ of Theorem~\ref{welldefined} one gets that our affine winding 
number $\win_p$ is defined for a vast collection of closed $M$ and connected 
components $\mathcal N$ of $C^{\infty}(N, M).$ Recall that for closed $M$ the 
winding number can not be defined as the particular case of the classical 
linking number invariant.  
 
\begin{thm}\label{compareclassical} 
Let $M$ be a manifold with $\partial M\neq \emptyset$ or such that it is the 
interior of a manifold with boundary. Let $\mathcal N$ be a connected 
component of $C^{\infty}(N, M)$ consisting of $f:N\to M$ with $f_*([N])=0\in 
H_*(M)$. $($This is the setup where we can define $\swin_p:\mathcal N\setminus 
\Sigma\to \Z$ as the linking number between $f(N)$ and $S^0$ consisting of $p$ 
and  a point in a boundary component.$)$ 
Then $\win_p-\swin_p:\mathcal N\setminus \Sigma\to \Z$ is a constant function. 
\end{thm} 
 
\pp 
Let $p^+$ denote the 0-dimensional singular cochain $1\cdot \gf$ where $\gf: 
\Delta^0 \to M, \gf(\Delta^0)=p$. We triangulate $N$ and regard a map $N \to M$ 
as a singular chain in $M$. Take $f\in \mathcal N\setminus \Sigma$ and recall 
that  $\swin_p(f)$ is defined as a particular case of the linking number: namely 
as the intersection number $\mathfrak S\bullet p^+$ where $\mathfrak S$ is a 
singular chain with $\partial \mathfrak S=f$. Let $F:N\times I\to M$ be a good 
homotopy 
between the preferred map $\eps\in \mathcal N\setminus\Sigma$ and $f.$ Take a 
singular chain $\overline {\mathfrak S}$ with boundary $\eps$ and consider 
the triangulation of $N\times I$ such that $F|_{N\times 0}$ is equal to 
$\partial \overline {\mathfrak S}$. Then $\mathfrak S:= \overline{\mathfrak S} + 
F$ is a singular chain with $\partial \mathfrak S=f$. Clearly  
$\mathfrak S\bullet p^+ -\overline{\mathfrak S}\bullet p^+=F\bullet p^+.$ By 
Remark~\ref{alter}, $F\bullet p^+=\Delta_{\win_p}(F).$ By Definition of 
$\swin_p$ 
we have 
$$
F\bullet p^+= {\mathfrak S}\bullet p^+-\overline{\mathfrak S}\bullet 
p^+=\swin_p(f)-\swin_p(\eps).
$$
Hence $\swin_p(f)-\swin_p(\eps)=\Delta_{\win_p}(F).$  
 
By Theorem~\ref{m=0}, $\A=\Z$ and $q=\id:\Z\to 
\Z$ for $M$ non-closed. By Definition of $\win_p,$ we have  
$$
\win_p(f)=q(\Delta_{\win_p}(F))=\Delta_{\win_p}(F)= 
\swin_p(f)-\swin_p(\eps). 
$$
Thus $\win_p-\swin_p:\mathcal N\setminus \Sigma\to \Z$ is the constant function 
 $-\swin_p(\eps).$ 
\qed 
 
\section{Some Applications}\label{apps} 
 
As a first application of our affine winding numbers, we have the following  
generalization of the Cauchy integral formula. 
 
\begin{thm}\label{Cauchy} 
Let $F^2$ be a $2$-dimensional surface equipped with a complex structure such  
that either $F$ is not closed or $F$ has genus bigger than one. Let $f$ be a  
meromorphic function on $F^2$ having poles $\{a_j\}_{j=1}^{k}$ and residues  
$\Res f(a_j), j=1,\dots, k$. Let $C_i:S^1\to F^2, i=1,2,$ be two homotopic 
smooth  
oriented $($not necessarily zero homologous$)$ curves not passing through any of 
the 
poles. Let  
$\win_{a_j}(C_i)\in \Z$, $i=1,2; j=1, \dots, k,$ be the affine winding numbers 
that are defined  for $p=a_j,$ since condition $0$ or condition $1$ of Theorem~$\ref{m=0}$ is satisfied. Then  
\begin{equation}\label{Cauchyidentity} 
\oint_{C_1} fdz=\oint_{C_2} fdz+2\pi i\Bigl (\sum_{j=1}^k \Res  
f(a_j)\bigl(\win_{a_j}(C_1)-\win_{a_j}(C_2)\bigr)\Bigr). 
\end{equation}   
\end{thm} 
 
\m 
This Theorem allows one to express the integral of a meromorphic function over  
a curve, through the integral of the function over some specified  
homotopic curve. The formulation of the classical Cauchy Theorem for $F^2=\C$ is 
obtained from  
Theorem~\ref{Cauchy} by taking $C_2$ to be a small curve lying far away from all  
the poles, so that $\oint_{C_2} fdz=0$. In this case by 
Theorem~\ref{compareclassical} 
$\bigl(\win_{a_j}(C_1)-\win_{a_j}(C_2)\bigr)$ coincides with the classical  
winding number $\swin_{a_j}(C_1).$  
 
Note that by Statement $2$ of Theorem~\ref{welldefined}, the term   
$\bigl(\win_{a_j}(C_1)-\win_{a_j}(C_2)\bigr)$ in~\eqref{Cauchyidentity} does not 
depend on the choice of the preferred map $\eps\in \mathcal N\setminus \Sigma$ 
used to define $\win_{a_j}=\win_{a_j, \eps}.$ 
 
\pp 
For non-closed $F$ the affine winding numbers $\win_{a_j}(C_i)$ are  
$\Z$-valued, since condition $0$ of Theorem~\ref{m=0} holds. 
For closed $F$ of genus bigger than one, the affine winding numbers 
$\win_{a_j}(C_i)$ are $\Z$-valued, since 
$b_1(S^1)+b_0(S^1)=2<b_1(F),$ and hence condition $1$ of \theoref{m=0} holds.

Similarly to the proof of the classical Cauchy Theorem, the proof of 
Theorem~\ref{Cauchy} boils down to local considerations. Namely,  
using statement $2$ of Theorem~\ref{welldefined} one shows that both parts of  
identity~\eqref{Cauchyidentity} change in the same way under an elementary  
homotopy of $C_1$ that involves one passage of $C_1$ through one of the poles. 
Since $C_1$ is homotopic to $C_2$ and the two sides of~\eqref{Cauchyidentity} 
are equal for $C_1=C_2,$ we get that identity~\eqref{Cauchyidentity} holds. 
\qed 
 
\subsection*{Applications of $\win_p$ to the study of wave front propagation.} 
Informally speaking,  
we assume that at a moment of  
time $T$ something happens at a submanifold of $M^m$ and the perturbation  
caused by this event starts to radiate from the submanifold in all the 
directions  
according to a propagation law. More accurately, we have a smooth map $W: 
N^{m-1}\times [T,  
\infty) \to M^m,$ where  $\Im W|_{N\times t}$ is thought of as the set of 
points that the perturbation has just reached at time $t.$ 
 
In fact, for wave fronts in geometric optics the map $W$ has  
special properties. For example $W\big|_{N\times t}, t\in [T,\infty),$ lifts to 
a Legendrian submanifold of the unit cotangent bundle of $M$, see 
Arnold~\cite{Arnold}.   
In this work we do not use any of these properties.  
 
We define the wave front $W(t): N \to M$ by setting  
$W(t)(n)=W(n,t), n\in N,$ and make an assumption that $W$ is {\em generic\/} 
i.e. $p$ is a regular value of $W$ and $(W(t))^{-1}(p)$ has at  
most one point for each $t\in[T,\infty)$.  
 
We would like to find an estimate from  
below on the number of  times $\pass(t_1, t_2)$ a wave front $W(t)$ on $M$ passed 
through  
the point $p$ between two moments of time $t_1$ and $t_2$ such that $p\not \in 
\Im W(t_i), i=1,2$.  
Moreover, we would like this estimate to be computable
from the shape of the pairs $(W(t_1), p)$ and $(W(t_2), p)$  
only, without any  
knowledge of $W$, topology of $M$, time  
moments $t_1, t_2$ etc. Clearly, we have  
$\pass(t_1,t_2)\ge|\Delta_{\win_p}(F)|,$ where  
$ 
F: N \times I \to M,\quad F(x,t)=W\left(x, (t_2-t_1)t+t_1\right). 
$ 
The difficulty is that we know $W(t_1), W(t_2),$ but we do not know $W(t)$ for 
$t_1<t<t_2$, and thus we do not know $F.$  
 
Luckily statement $2$ of Theorem~\ref{welldefined} says that
$ 
q(\Delta_{\win_p}(F))=\win_p(W(t_2))-\win_p(W(t_1)) 
$ and that we  
can take {\it any} good homotopy $G$ between $W(t_1)$ and $W(t_2)$ and compute 
$\win_p(W(t_2))-\win_p(W(t_1))$ as $q(\Delta_{\win_p}(G))$. In particular, if 
$M, N, \mathcal N$  
are as in~\theoref{m=0}, so that $\A=\Z$ and $q=\id:\Z\to \Z,$  then
$\pass(t_1, t_2)\geq |\win_p(W(t_2))-\win_p(W(t_1))|.$ Thus in this case we  
can estimate from below $\pass(t_1,t_2)$ from the pictures of the front at times 
$t_1$ and  
$t_2$.   
 
\begin{figure}[htbp]  
\begin{center}  
\epsfxsize\hsize\advance\epsfxsize -0.5cm  
\epsfxsize 12 cm  
\hepsffile{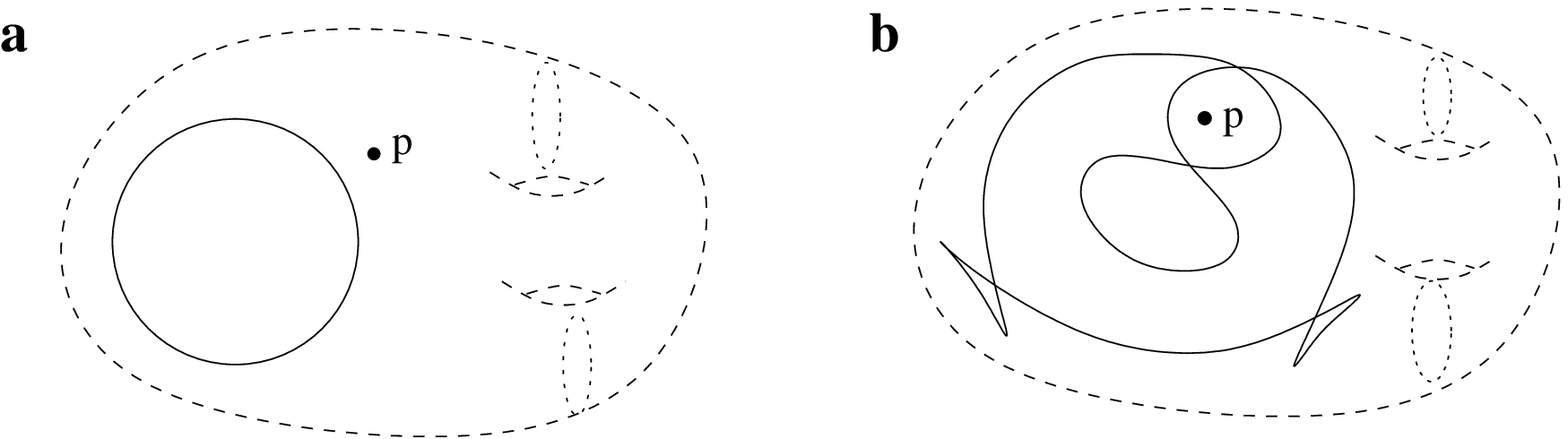}  
\end{center}  
\caption{}\label{examplewinding.fig}  
\end{figure}

\begin{ex}\label{example}  
Assume that at times $t_1$ and $t_2$ the wave front is contained in a chart of 
$M.$ 
Assume moreover that at time $t_1$ the picture of the  
wave front was the one shown in Figure~\ref  
{examplewinding.fig}a and later at $t_2$  
it developed into the shape shown in Figure~\ref  
{examplewinding.fig}b.  
(The Figure~\ref{examplewinding.fig}b depicts a sphere  
that can be obtained  
from the trivially embedded sphere  
by passing two times through the point $p$ and by creation of some  
singularities on the part of the front away from $p$.)  
 
A straightforward  
calculation gives 
$  
\win_p(W(t_2))-\win_p(W(t_1))=q(\pm 2)\in \A,$ where the sign depends on
the front orientation which is not shown in the Figure.  
 
Assume that  $M$ is not a rational homology sphere, then 
condition $3$ of Theorem~\ref{m=0} is satisfied. Hence $\A=\Z$ and $|q(\pm 
2)|=2=|\win_p(W(t_2))-\win_p(W(t_1))|.$  We conclude that every generic $W$ that changes $W(t_1)$ to  $W(t_2)$ involves at least two 
passages through $p$.  
 
If $W$ is not generic, it could happen that two branches of 
the front pass through $p$ simultaneously. However for non-generic $W$ we still can conclude that the front passed through $p$ at least 
once between the two time moments. 
\end{ex} 
 
\section{Another generalization of the winding number and  
the $\aawin$ invariant.}\label{versiya}  
 
This section deals with the case where the point, around which we compute the 
winding number, is not stationary but rather moves in $M.$ We will define the 
affine winding number  
invariant similarly to how we did it before, but the indeterminacy  
subgroup $\bA$ will increase.  
 
As before we fix a connected component $\mathcal N$ of $C^{\infty}(N, M).$  
We put the {\it discriminant\/} $\wt \Sigma\subset \mathcal N\times M$ to be the 
set of pairs $(f, x)\in \mathcal N\times M$ such that $x\in \Im f.$ Note that we 
do not include into the discriminant the pairs $(f, x)$ such that $f$ is 
singular in the common sense but $x\not \in \Im f.$

Let $F:N \times  I \to M$ be a smooth map and let $\gamma:I\to M$ be a smooth 
path such that $\gamma(t)\cap F(N\times t)=\emptyset, t=0,1.$ By a 
$C^{\infty}$-small perturbation we may assume that $F$ and $\gamma$ are 
transverse. 
 
We define $\Phi: N\times I \to M\times I, \Phi(n,t)=(F(x, t),t)$ and $\Gamma: I 
\to M\times I, \Gamma (t)=(\gamma(t),t)$. Since $F$ and $\gamma$ are transverse, 
we see that $\Phi$ and $\Gamma$ are also transverse. 
We equip $I=[0,1]$ with the orientation from $0$ to $1$ and we equip $N  
\times I$ and $M\times I $ with the product orientations.  
 
\begin{defin}\label{newDelta}  
We define $\Delta_{\aawin}(F,\gamma)\in \Z$ as the intersection number 
$\Phi\bullet \Gamma$ of $\Phi$ and $\Gamma$ with respect to the orientations 
described above.  
 
Note that if $\gamma$ is a constant path, so that $\Im \gamma$ is just one point 
$p\in M,$ then $\Delta_{\aawin}(F,\gamma)=\Delta_{\win_p}(F),$ see 
Remark~\ref{alter}. 
\end{defin} 
 
\begin{defin}[$\aawin(f, x)$-invariant]\label{wtwindefin} 
Fix a pair $(\eps, \rho)=(\eps_{\mathcal N}, \rho)\in\mathcal N\times M\setminus 
\wt \Sigma$ that should be thought of as a preferred map and a preferred point. 
 
Let $\bA\subset \Z$ be the subgroup from \defref{A},  
and let $\bB\subset \Z$ be the subgroup  consisting of numbers  
$\eps_*([N])\bullet a\in \Z,$ where $\bullet$ denotes the intersection pairing 
and $a$  
runs over $H_1(M)$. Put $\B=\Z/(\bA+\bB)$ and put $\wt q:\Z\to \B$ to be the 
quotient homomorphism. Clearly $\B$ does not depend on the choice of $(\eps, 
\rho).$

Given $(f, x)\in \mathcal N\times M\setminus \wt \Sigma,$ we define the {\it 
affine winding number\/} $\aawin(f, x)=\aawin_{\eps, \rho}(f, x)\in \B$ by 
setting $\aawin(f, x)=\wt q(\Delta_{\aawin}(F, \gamma))\in \B,$  
where $F:N\times I \to M$ is a good homotopy between $\eps=F|_{N\times 0}$ and 
$f=F|_{N\times 1}$ 
and $\gamma$ is a smooth path from $\rho$ to $x$ that is transverse to $F.$  
\end{defin} 
 
\begin{thm}\label{wtwinwelldefined}~\begin{description} 
\item[1] 
The number $\aawin(f,x)=\aawin_{\eps, \rho}(f, x)\in \B$ does not depend on the 
choices of the good homotopy $F$ and of the path $\gamma$ in its definition. 
Thus the function
$$
\aawin:\mathcal N\times M\setminus \Sigma \to \B
$$ 
is well defined.  
\item[2] Take $(f_0, x_0), (f_1, x_1)\in \mathcal N\times M\setminus \wt 
\Sigma,$ a good homotopy $H:N\times I\to M$ between $f_0=H|_{N\times 0}$ and 
$f_1=H|_{N\times 1},$ 
and a path $\beta$ from $x_0$ to $x_1$ that is transverse to $H.$   
Then 
$$
\aawin(f_1, x_1)-\aawin(f_0, x_0)=\wt q(\Delta_{\aawin}(H, \beta)).
$$ 
In 
particular, $\aawin:\mathcal N\times M\setminus \wt \Sigma \to \B$ is a locally 
constant function. 
\item[3] If we replace $(\eps, \rho)$ by some $(\eps', \rho')\in \mathcal 
N\times M\setminus \wt \Sigma,$ then  
$$
\aawin_{\eps', \rho'}-\aawin_{\eps, \rho}:\mathcal N\times M\setminus \wt \Sigma\to \B
$$ 
is the constant function $\aawin_{\eps', \rho'}(\eps, \rho).$ Hence the 
functions 
$$
\aawin_{\eps', \rho'}, \aawin_{\eps, \rho}:\mathcal N\times M\setminus \wt 
\Sigma\to \B
$$ 
are equal up to an additive constant. 
\end{description} 
\end{thm} 
 
\pp {\em Let us prove statement $1$.\/} Choose another good homotopy $F'$ of 
$\eps$ to $f$ and another path $\gamma'$ from $\rho$ to $x$ that is transverse 
to $F'.$ Let $\Phi':N\times I\to M\times I$ and $\Gamma':I\to M\times I$ be the 
corresponding maps. 
Clearly it suffices to show that $\Phi\bullet \Gamma-\Phi'\bullet \Gamma'\in 
\bA+\bB.$

Consider a non-decreasing smooth function $\gf: [0,1]\to [0,1]$ such that 
$\gf(t)=0,$  
for $t\in[0,1/2],$ and $\gf(1)=1$. Define the good homotopies $\ov F, \ov F'$ 
via 
$\ov F(x,t)=F(x,\gf(t)), \ov F'(x,t)=F'(x,\gf(t))$ and put $\ov \Phi, \ov 
\Phi':N\times I\to M\times I$ to be the corresponding maps. Define $\ov \gamma, 
\ov \gamma':I\to M$ via   
$\ov \gamma(t)=\gamma(1-\gf(1-t)), \ov \gamma'(t)=\gamma'(1-\gf(1-t))$ and put 
$\ov \Gamma, \ov \Gamma':I\to M\times I$ to be the corresponding maps. 
 
Since $\Phi$ is homotopic to $\ov \Phi$ and $\Gamma$ is homotopic to $\ov 
\Gamma$ modulo boundary, we have $\Phi\bullet \Gamma=\ov \Phi \bullet \ov 
\Gamma,$ where $\bullet$ is the intersection pairing of the corresponding 
relative homology classes modulo $\partial (M\times I).$ 
Similarly $\Phi'\bullet \Gamma'=\ov \Phi' \bullet \ov \Gamma'.$  
Thus it suffices to show that $\ov \Phi\bullet \ov \Gamma-\ov \Phi'\bullet \ov 
\Gamma'\in \bA+\bB.$ 
 
Since the maps $\ov \Phi, \ov \Phi',\ov \Gamma, \ov \Gamma'$ preserve the 
$I$-coordinate, we have  
that  
\begin{equation}\label{long} 
\begin{split} 
\ov \Phi\bullet \ov \Gamma-\ov \Phi'\bullet \ov \Gamma'= 
\bigl(\ov \Phi|_{N\times [0,\frac{1}{2}]}\bullet \ov 
\Gamma|_{[0,\frac{1}{2}]}+\ov \Phi|_{N\times [\frac{1}{2}, 1]}\bullet \ov 
\Gamma|_{[\frac{1}{2},1]}\bigr)-\\\bigl(\ov \Phi'|_{N\times 
[0,\frac{1}{2}]}\bullet \ov \Gamma'|_{[0,\frac{1}{2}]}+\ov \Phi'|_{N\times 
[\frac{1}{2},1]}\bullet \ov \Gamma'|_{[\frac{1}{2},1]}\bigr)= 
\mathfrak A+\mathfrak B, 
\end{split} 
\end{equation} 
where  $\mathfrak A=\ov \Phi|_{N\times [\frac{1}{2}, 1]}\bullet \ov 
\Gamma|_{[\frac{1}{2},1]}-\ov \Phi'|_{N\times [\frac{1}{2},1]}\bullet \ov 
\Gamma'|_{[\frac{1}{2},1]}$ and  $\mathfrak B=\ov \Phi|_{N\times 
[0,\frac{1}{2}]}\bullet \ov \Gamma|_{[0,\frac{1}{2}]}-\ov \Phi'|_{N\times 
[0,\frac{1}{2}]}\bullet \ov \Gamma'|_{[0,\frac{1}{2}]}.$ 
 
By our choice of $\gf,$ the restrictions $\ov \gamma|_{[\frac{1}{2}, 1]}, \ov 
\gamma'|_{[\frac{1}{2}, 1]}$ are constant paths from $x$ to $x$ and $\ov F 
|_{N\times [\frac{1}{2}, 1]},  
\ov F'|_{N\times [\frac{1}{2}, 1]}$ are smooth homotopies of $\eps$ to $f.$ 
Define the good homotopies $\wt F, \wt F':N\times I\to M$ of $\eps$ to $f$ via 
$\wt F(n, t)=\ov F(n, \frac{1}{2}t), \wt F'(n, t)=\ov F'(n, \frac{1}{2}t),$ for 
$n\in N, t\in I.$ 
By Remark~\ref{alter} and the observation in Definition~\ref{newDelta} we get 
that $\mathfrak A=\Delta_{\win_{x}}(\wt F)-\Delta_{\win_{x}}(\wt F').$ By 
statement $2$ of Theorem~\ref{welldefined} we have that 
\begin{equation}
\begin{aligned}
0&=\win_{x}(f)-\win_{x}(f)=q(\Delta_{\win_{x}}(\wt F))- 
q(\Delta_{\win_{x}}(\wt F'))\\
&=q(\mathfrak A)\in \A=\Z/\bA. 
\end{aligned}
\end{equation}
Hence $\mathfrak A\in  \bA$. 
 
By our choice of $\gf,$ the restrictions $\ov F|_{N\times [0,\frac{1}{2}]}$ and 
$\ov F'|_{N\times [0,\frac{1}{2}]}$ do not move $N$ on $M$ and $\ov 
\gamma|_{[0,\frac{1}{2}]}, \ov \gamma'|_{[0,\frac{1}{2}]}$ are paths from $\rho$ 
to $x.$ Thus $\mathfrak B$ equals to the intersection number of $\eps$ and of 
the closed loop   $\ov \gamma|_{[0,\frac{1}{2}]}\left(\ov 
\gamma'|_{[0,\frac{1}{2}]}\right)^{-1}$. Hence $\mathfrak B\in \bB.$ 
 
So, $\mathfrak A +\mathfrak B \subset \bA+\bB$. Now identity~\eqref{long} implies that $\ov \Phi\bullet 
\ov \Gamma-\ov \Phi'\bullet \ov \Gamma' \in \bA+\bB$ and 
we proved statement $1$ of the Theorem.  

\m 
The proofs of statements $2$ and $3$ are similar to the proofs of the 
corresponding statements of Theorem~\ref{welldefined} and therefore are omitted. 
\qed

The following Theorem says that $\aawin$ is a $\Z$-valued invariant 
for many manifolds $M$ and connected components $\mathcal N$ of $C^{\infty}(N, 
M).$  
 
\begin{thm}\label{wtwinZ}  
The equality $\B=\Z$ holds if $\mathcal N$ consists of maps $g$ such that 
$g_*([N])\in H_*(M)$ 
is a finite order element and at least one of the conditions {\bf 0,1,2,3,4,5\/} 
of  
Theorem~$\ref{m=0}$ is satisfied.  
\end{thm}  
 
\pp By Theorem~\ref{m=0}  $\bA=0\subset \Z.$ Since $\eps_*([N])\in H_*(M)$ is an 
element of finite order, we have $\eps_*([N])\bullet a=0$ for every $a\in H_1(M).$ 
Hence $\bB=0\subset \Z.$ Thus $\B=\Z/(\bA+\bB)=\Z.$\qed

\begin{rem}[Comparison of $\aawin$ to the winding number that is defined as a 
particular case of the linking number]\label{wtwincompareclassical} 
Let $M$ be a manifold with $\partial M\neq\emptyset$ or that is the interior of 
a manifold with boundary. Let $\mathcal N$ be a connected component of 
$C^{\infty}(N, M)$ consisting of $f:N\to M$ with $f_*([N])=0\in H_*(M).$ This is 
the setup where we can define $\swin:\mathcal N\times M\setminus \wt \Sigma\to 
\Z$ as the linking number between $f$ and the map of $S^0$ sending one point of $S^0$ to $p$ and the other into a chosen boundary component. 
Then $\aawin-\swin:\mathcal N\times M\setminus \wt \Sigma\to \Z$ is a constant 
function and hence the two invariant are equal up to an additive constant. 
The proof of this fact is similar to the proof of Theorem~\ref{compareclassical} 
and therefore is omitted. 
\end{rem}

\begin{rem}[$\aawin$ and Viro's winding number] In~\cite{Viro} Viro introduced 
generalizations of winding numbers to the case of zero homologous immersed 
curves on a closed surface $F^2$ with $\chi(F^2)\neq 0$. Viro's winding numbers 
are $\Q$-valued and under regular homotopy they behave in the same way as the 
classical winding numbers. However unlike our affine winding numbers, Viro's 
winding number around $p\in F^2$ changes under a non-regular homotopy of a curve 
that does not pass through $p$. Hence Viro's winding number does not give rise 
to a locally constant function on $\mathcal N\times F^2\setminus \wt \Sigma.$ 
\end{rem}

\begin{rem}[$\aawin$ and wave fronts]  
The invariant $\aawin$ allows one to estimate from below the number of times a 
front 
on $M$ passed through a point that was continuously moving in $M.$ Assume that 
the pair: the trajectory $\alpha:[T, +\infty)\to M$ of the point and the front 
propagation  $W:N^{m-1}\times [T, +\infty)\to M^m$ is generic so that $\alpha$ 
and $W$ are transverse. Choose $t_1, t_2>T$ such that $\alpha(t_i)\not \in \Im 
W|_{N\times t_i}, i=1,2,$ and define  
$F: N \times I \to M,\quad F(x,t)=W\left(x, (t_2-t_1)t+t_1\right), \quad 
\gamma:I\to M,\quad  \gamma(t)=\alpha 
((t_2-t_1)t+t_1).$ Clearly 
$|\Delta_{\aawin}(F, \gamma)|$ estimates from below the number of times the 
front passes through the moving point between times $t_1$ and $t_2.$  
 
Put $f_i=W|_{N\times t_i}:N\to M, i=1,2.$ 
If $M$ and $\mathcal N\ni W|_{N\times T}$ are such that $\aawin$ is a 
$\Z$-valued invariant, then $\Delta_{\aawin}(F, \gamma)=\aawin(f_2, 
\gamma(t_2))-\aawin(f_1, \gamma(t_1)).$ By statement  
$2$ of Theorem~\ref{wtwinwelldefined} the last quantity can be computed using 
any good homotopy $H$ of $f_1$ to $f_2$ and a path $\beta$ from $\gamma(t_1)$ to 
$\gamma(t_2)$ that is transverse to $H.$ 
This allows us to estimate from below the number of times a front passed through 
an observable point moving on $M$ between times $t_1$ and  $t_2.$ This 
estimation can be done from the snapshots of the front and the point at the two 
time moments without the knowledge of the front and point movements.  
 
For example, assume that at times $t_1$ and $t_2$ the front and the observable 
point were located in a chart of $M$ and were as depicted in 
Figure~\ref{examplewinding.fig}. Assume that $M$ is not a $\Q$-homology sphere, 
so that $\B=\Z$ by Theorem~\ref{wtwinZ}. By the above discussion we get that for 
generic $(W, \alpha)$ the front passed through an observable moving point at 
least $2$ times between times $t_1$ and $t_2.$ 
{\em  This conclusion can be made without  
the knowledge of $W$ and of the trajectory $\alpha$ of the observable point.\/} 
Similarly to the case of a stationary point $p$ discussed in 
Example~\ref{example},  for non-generic pairs $(W, \alpha)$ we still can 
conclude that the front passed through the moving point at least once between 
the two time moments. 
\end{rem}

{\bf Acknowledgment:} The first author was partially supported by the Walter and 
Constance Burke Research Initiation Award. 
The second author was partially supported by NSF, grant 0406311, USA and by 
MCyT, projects BFM 2002-00788 and BFM2003-02068/MATE, Spain. His visit to 
Dartmouth College was supported by the funds donated by Edward Shapiro to the 
Mathematics Department of Dartmouth College.


\begin{thebibliography}{99999}  
\bibitem{AMR1} 
J.~E.~Andersen, J.~Mattes, N.~Reshetikhin, {\em  
Quantization of the algebra of 
chord diagrams,\/} Math. Proc. Cambridge Philos. Soc., Vol. {\bf 124} no. 3 
(1998), pp. 451--467  
 
\bibitem{AMR2} 
J.~E.~Andersen, J.~Mattes, N.~Reshetikhin, 
{\it The Poisson structure on the moduli space of flat connections and chord 
diagrams.\/} Topology {\bf 35} (1996), no. 4, 1069--1083. 
 
 
 
\bibitem{Arnold}  
V.I.~Arnold: Invariants and perestroikas of wave fronts on  
the plane, Singularities of smooth maps with additional  
structures,  
{\em Proc. Steklov Inst. Math.}, Vol. {\bf 209} (1995),  
pp. 11--64.  
 
 
\bibitem{Tchernov}  
V.~Chernov (Tchernov): {\em Shadows of wave fronts and Arnold-Bennequin type 
invariants of fronts on surfaces and orbifolds.\/} Differential and symplectic 
topology of knots and curves, 153--184, Amer.~Math.~Soc.~Transl.~Ser. 2, {\bf 
190}, Amer.~Math.~Soc., Providence, RI (1999) 
  
             
  
 
\bibitem{ChernovPoisson} 
V.~Chernov (Tchernov), {\it Graded Poisson algebras on bordism groups of 
garlands and their applications,\/} available as a preprint math.GT/0608153 at 
the http://www.arXiv.org (2006) 
 
 
 
 
\bibitem{CR} 
V.~Chernov (Tchernov) and Yu. B. Rudyak: {\em Toward a  
General Theory of  
Affine Linking Invariants,\/} Geometry and Topology, Vol. 9 (2005) Paper no. 42,  
pages 1881--1913; http://www.maths.warwick.ac.uk/gt/GTVol9/paper42.abs.html  
 
 
\bibitem{ChernovRudyakGarlands}  
V.~Chernov (Tchernov) and Yu.~B.~Rudyak, {\it Algebraic Structures on  
Generalized Strings,\/} available as a preprint  
math.GT/0306140 at the http://www.arXiv.org (2003) 
 
 
 
\bibitem{DoCarmo}  
M.~P.~do Carmo: {\em Riemannian Geometry}, Birkhauser,  
Boston, 1992  
 
\bibitem{Goldman}  
W.~Goldman: {\it Invariant functions on Lie groups and  
Hamiltonian flows of surface group representations.\/}  
Invent. Math. {\bf 85} (1986), no. 2, 263--302.  
 
 
  
\bibitem{Goryunov}  
V.~Goryunov: {\em Local invariants of mappings of surfaces into three-space.\/} 
The Arnold-Gelfand mathematical seminars, pp.~223--255, Birkhäuser Boston, 
Boston, MA (1997)  
 
 
\bibitem{Kaiserbook} 
U.~Kaiser, {\it Link theory in manifolds.\/} Lecture Notes  
in Mathematics, 1669.  
Springer-Verlag, Berlin, (1997)  
 
 
\bibitem{MikhalkinandPolyak} 
G.~Mikhalkin and M.~Polyak: {\em Whitney formula in higher dimensions.}  
J. Differential Geom. 44 (1996), no. 3, pp.~583--594.  
 
 
\bibitem{Polyak} 
M.~Polyak: {\em Shadows of Legendrian links and $J\sp +$-theory of curves.\/}  
Singularities (Oberwolfach, 1996), pp.~435--458,  
Progr. Math., {\bf 162},  
Birkhäuser, Basel (1998)  
 
 
\bibitem{Preissman} 
A.~Preissman: {\em Quelques propri\'et\'es globales des espaces de  
Riemann}, 
Comment. Math. Helv. {\bf 15} (1943), pp.~175--216. 
 
\bibitem{Rudyak}  
Yu.~B.~Rudyak: {\em On Thom Spectra, Orientability, and  
Cobordism}, Springer,  
Berlin Heidelberg New York (1998)  
 
\bibitem{ShumakovichStrangeness} 
A.~Shumakovich: 
{\em Explicit formulas for the strangeness of plane curves.\/} (Russian. Russian 
summary)  
Algebra i Analiz {\bf 7} (1995), no. 3, 165--199; translation in St. Petersburg 
Math.~J.~{\bf 7} (1996), no. 3, pp.~445-472 
 
\bibitem{ShumakovichVassiliev} 
A.~Shumakovitch: {\em Shadow formula for the Vassiliev invariant of degree 
two.\/} Topology {\bf 36} (1997), no. 2, 449--469 
  
             
\bibitem{Turaevshadow}  
V.~Turaev: {\em Shadow links and face models of statistical mechanics.\/} 
J.~Differential Geom.~{\bf 36} (1992), no. 1, pp.~35--74 
 
\bibitem{Turaevskein}  
V.~Turaev: {\it Skein quantization of Poisson algebras of  
loops on surfaces.\/} Ann. Sci. Ecole Norm. Sup.  
(4) {\bf 24} (1991), no. 6, pp.~635--704.  
 
 
 
 
\bibitem{Viro} 
O.~Viro: {\em Generic immersions of the circle to surfaces and the complex 
topology of real algebraic curves.\/} Topology of real algebraic varieties and 
related topics, pp.~231--252, Amer.~Math.~Soc.~Transl.~Ser.~2~{\bf 173,} Amer. 
Math. Soc., Providence, RI (1996)   
 
 
\bibitem{Whitney}  
H.~Whitney: {\em On regular closed curves in the plane.\/} 
Compos. Math. 4, 276-284 (1937)  
 
 
\end{thebibliography}
\end{document}